# What is Math Really?

Anthony Rizzi

*Institute for Advanced Physics, Baton Rouge, LA 70895*

**Abstract:** Modern mathematics is known for its rigorous proofs and tight analysis. Math is the paradigm of objectivity for most. We identify the source of that objectivity as our knowledge of the physical world given through our senses. We show in detail, for the core of modern mathematics, how modern mathematical formalism encapsulates deep realities about extension into a system of symbols and axiomatic rules. In particular, we proceed from the foundations in our senses to the natural numbers through integers, rational numbers, and real numbers, including introducing the concept of a field. An appendix shows how the formalism of complex numbers arises. [1]

## I. Introduction

We are all taught math from a young age, but, practically none of us is ever told what math is truly all about. We learned how to do it, what worked and what did not, but not usually why, and seldom where any of it came from. As we moved from arithmetic to fractions to algebra, any thoughts we had of where math came from became more and more hazy. Sets, axioms and proofs begin to take center stage. Classrooms around the country in middle and high school years ring with "What will I ever use this for?" and "What does this have to do with life?" To most, math is now simply a tool for certain technical professions, not something of any real interest or relevance. Such professions are considered to be extremely important, but very hard and best reserved for those trained technicians who can do such things without hurting their brains too much. (Or, there is always a sneaking suspicion, maybe they are just more willing to do so?) They are thought to be like the acrobats at the circus who jump through hoops of fire and stand on thin wires or like the lineman fixing the high voltage equipment on telephone polls high in the air, but without any of the excitement of these activities.

By contrast, someone like myself, who grows up and does enter such a mathematical profession (I am a physicist), feels comfortable saying, as the renowned physicist Eugene Wigner did, how amazing is "the unreasonable effectiveness" of mathematics in revealing a deep understanding about the world. Those who take up mathematics for a profession will often think of math as existing out there in some realm of unchanging things, which Plato so aptly expressed. Other mathematicians will see math as it is treated in formal proofs that dominate so much of modern mathematics, that is, as simply the necessary consequences of a system of axioms chosen freely. Bertrand Russell took this to its ultimate conclusion when he said: "Mathematics [is the science] in which we never know what we are talking about, nor whether what we are saying is true."[2] How indeed does this seemingly otherworldly, "not part of real life" thing in fact have so much to do with the real world?

It sure would be nice to know what math is, for it is always nice to know what one is doing, especially for those of us who do so much of it. Though we have largely forgotten, Aristotle [1] answered the question "What is math?" in its most rudimentary form over 300 years before Christ. His answer, which we will come to briefly, is profound and simple, as such primary things are. Hence, it is critically important to teach children the answer at its most basic level from a young age and to reinforce it and deepen it as they grow in their understanding. Truly understanding where math comes from and what it is about will certainly mean less frustration over and more attention to the subject. As a result, both those who go far and those who do not will go further than they do now in the queen of the sciences. In this paper, we will start with kindergarten math and work through high school math bringing out how we go from the concrete reality around us to the axioms of mathematics.

---

[1] This paper was written in 2006. Appendix was written in 2018.
[2] B. Russell, *Mysticism and Logic and Other Essays* (Longmans, Green & Co., London, 1919), pg 75.

There are two ways we know: through the senses and through the mind; the second is dependent on the first. We first sense something, then we abstract, or pull something out of it, with our mind. In math, we consider the abstract. We discuss lines in general, not this particular line. Indeed, even if we try to focus on a particular line, say of a line of a given length, we realize that it is also abstract. We have left out-- abstracted away-- the height and width we see in things to consider only their length.

Aristotle pointed out that mathematics is the study of *quantity*. Quantity is what we get when we take what we see and leave all of its properties behind but its *extension*. To be concrete, pick an object, freeze time, take away the environment around the object, take away its abilities to act, its color, its hardness and its smell; strip all away but extension and what limits the extension.[3] So, for example, if I look at a round pie and abstract all away but extension, I am left with a circular disk. We call the boundary of the disk a circle. This is basic geometry. It is what we do in Kindergarten and before. You remember or you have seen your child or another child doing it: "Identify the circle in this group of figures;" "Identify the square." It is the very important beginning stage. But, counting comes quickly also. That is more what we think of as our first real math. Counting gives the "natural numbers."

In this way, we establish the two major subdivisions of math: geometry and arithmetic, which are roughly speaking, respectively, the study of extended shapes and the study of number. Let us examine each more closely.

## II. Arithmetic and Geometry

Arithmetic is more than counting. To see what more, consider the sensorial image, say, of a man. We know the man is a whole with parts; he is not ten fingers and ten toes, etc., but a man with ten fingers and ten toes. In *thinking* about the man, as opposed to simply considering his sensorial appearance as such, we must first come to see that he is something, a whole something, before we can ask about the properties of this something, such as his parts. In short, we first realize that there is a whole, then we understand the parts of the whole. For instance, we see he has two eyes. Each eye is a kind of unit, which I distinguish by leaving behind the accidental differences (e.g., one is darker brown or sees better than the other) and only keep the essence of being an eye (i.e. that it can see). In this way, I can see there is such a thing as two, that is, a whole which has two parts. This, however, is not yet sufficient to get the number "two" in its primal sense. To do this, we can imagine a whole which has two actual simple parts, i.e. actual parts which themselves do not have parts. Note that, in picturing such a thing, our

---

[3] For more discussion, see A. Rizzi, *A Kid's Introduction to Physics* (IAP press, Baton Rouge, 2012)

imagination is using realities we have already seen, not creating anything fundamentally new.

Let us pick a concrete simple image to illustrate the point. Imagine a circle with two parts. The circle is the whole and the diameter line defines the parts. From this image, you can abstract, pull out, the notion of two: "a whole composed of two units that do not share a common boundary." Your mind leaves out the extension and shape (and thus the boundaries) and considers only the units and the whole. Hence, the "do not share a common boundary" recalls that we have left behind the fact that the parts are touching and are only considering the whole with two parts as such. In this way, we get the abstract concept of two. It is now simply two, not two halves of a circular continuum, so it can, for example, be applied to two parts that divide the area of the circle in any way, or two parts of a sphere or whatever.

After getting this general idea in your mind, you can go back to the particular circle in your imagination, and divide the circle further by putting a radial line through one of the two parts, yielding a circle with three parts. These parts are not equal in terms of size, but we can leave out the size and think of the part as simply a primary unit: "a simple part that shares no boundary with another part." Each primary unit is identical in that each is an undivided part of the whole.

In this way, we see the nature of number in the most *primary* sense of the word. ***A number is a whole whose parts do not share a common boundary***; we also can call this discrete quantity. Note this is different from what mathematicians call the natural numbers, for the natural numbers include one. Though it is not often discussed and is little known, this concept of number, *not* the mathematician's "natural numbers", is, as we have seen, the one we come to first. Ancient authors such as Aristotle [1] were familiar with these distinctions as also were the innovators, such as Descartes [2], that put modern math on its successful highly axiomatic and symbolized path. We will call the above primal concept of number the ***fundamental or primary numbers, $N_F$***.[4]

We can make a simple jump to *analogically* generalize this primal concept of number. The primal concept of number just discussed is really abstracted from physical things and thus, in that way, belongs to them. For example, a given man, a substantial whole, really has so many (a number) actual parts, though we do not know how many. By contrast, if we just focus, as earlier, on his two eyes, we are considering just two parts out of the many parts of the full whole and each part is a unit only in some secondary sense. In short, because these parts are not themselves part-less and because they are not the sum total of the parts of the whole or even any continuous whole,[5] our use of the number two is here only analogical, i.e., only like the primal meaning of "two" in a limited way. Lumping these two meanings into a single concept "number" means the new concept "number" now is not simply from the physical world because it has an element of my mental work of lumping things together. These two concepts of number are really distinct

---

[4] We can reference the species in this genus of number (*fundamental concept)* symbolically as: $N_F = \{2, 3, 4...\}$.

[5] The type of analogical whole indicated in our original usage of "two eyes" is a very particular type of whole related to the powers of the man; in particular, the two eyes refer to that part of the qualitative power of seeing which consists of the organs with which the man sees. However, here, to stay in the realm of geometry, we would want to distinguish the eyes by their relative position in the body.

but we treat them as one for convenience of focusing on what is the same in the two concepts.

Another analogical generalization that goes further in ignoring both the distinctions between part-less parts and the parts being parts of a true whole, is to consider two things, say two apples, because, for example, they are in the same vicinity, as also being instances of the number "2." Note that such analogical generalizations are very powerful and useful but they, of course, do not make the distinctions that they ignore not real.

Keeping in mind this caveat and the dangers that attend forgetting it, we nonetheless can and should go still further; we can *analogically* generalize the concept of number so that it also applies to fundamentally different types of wholes, different than the property of geometric extension from which we originally pulled it. For example, I can obviously count a "number" of qualities of a given substance. I can, for instance, note that a copper bell is reddish-brown, melodious, hard, and cold; I count four properties. Yet, here I am using "number" in a different way because here we have a different type of "*whole with simple parts*," for the reddish-brown and hardness are there together on top of each other, not separated like two toes or two eyes.[6] Thus, we see counting properties of a bell is different than counting eyes, though we say we have a certain *number* of eyes and a certain *number* of properties, using the same word. Again, this applicability in many analogically similar types of cases makes the concept of number a very powerful and useful one.

*Geometry*, in the primary sense we consider here, is less abstract, leaves less behind; it studies wholes whose parts *do* share a common boundary such as the divided circle. Notice again that once we talk about 3, we have lost all idea of the boundary that exists in any real thing we see; we have left it behind in our notion. Whereas *arithmetic* studies discrete quantity (number), geometry is the study of what we can call continuous quantity.[7]

# III. Number

In summary, starting from full physical reality that we see, we leave behind all but the continuous quantity, i.e. the extension of the object. From there, we move to discrete quantity, which is the most abstract thing we can come to by such a process of abstraction. We start with a physical whole in our imagination (or in front of us). We see we can divide that whole, giving a new whole which has two parts, from which we get two. We can divide one of the parts of the physical whole again and get a whole with three parts. It is a little like embryonic cell division. From this, we see that every number is a whole that has been completed (or terminated) by a last unit (part). The last unit is the one that separates it from the number that comes before it, so that *three* has been completed by the last unit that makes it three which is not present in *two*. This last unit is arbitrary once we conceive the numbers in our mind (as opposed to imagining pictures of them), for the parts of the number *three* are indistinguishable; it is only when we compare it to *two* that we see one of the parts of *three* must be assigned as the last to complete the whole and make it *three* rather than *two*.

---

[6] The astute reader will notice that there is also a difference between the whole used for the sliced pie and the whole used for the eye as part of the man in the previous section. There, eye was used as an operation of the man (see footnote 4), rather than simply a geometric part of a larger shape as we could have had by using "eye" to mean only the part in the unique spatial location in the body (the sockets of the skull). Again, the power of the simple concept of number to be extended to many (analogically) different cases is so great, we sometimes do not even notice how different the cases are.

[7] Later, using principles we have pulled from what we see, we can piece together axioms that represent abstract geometry that we have not seen. This is part of the power of axiomatizing the principles, but to ensure consistency, we have to have a path back to the real.

Strange though it may sound, one is not a number in the proper sense, but is the unit by which we measure such wholes (numbers), *for it has no parts*. It sounds strange because the present use of the word number does not explicitly draw the distinction between the unit and the wholes composed of units. The fundamental unit is simple and not divided and is that by which the complex wholes are measured.[8] This is a real distinction, even when we do not explicitly make it. In leaving it behind, we leave part of the real behind. Even so, *implicit* distinctions between "one" (that *by which*) and the fundamental numbers (*that which* is measured) remain that all will recognize: for example, in the definition of primes, one is excluded, and is excluded also in axioms such as that of the multiplicative identity element and multiplicative inverse.

We have seen that the concept of number comes from the capacity of matter to be divided, i.e. cut into two with no limit; it comes from the very continuity that is in quantity. Again, number is abstract, meaning it leaves behind the continuity (the connectedness of one part to the next) that we see in real things.

Now, as interesting and profound as this is, it is not at all the way modern mathematics is done. Indeed, modern math is not explicitly thought about in this way. With this, we can see why the need arises for an article such as this to answer the question: "What is math?". The need arises because modern mathematical thinking does not focus on, indeed neglects, the very understanding (such as the generic introduction to geometry and number given above) that is the core of what math is. Still, modern mathematics *implicitly* takes its starting point in the ground given above (and continues to implicitly import aspects of extension given through our senses as it advances). We must start somewhere. Modern mathematics quickly moves to axiomatization of the principles that we have pulled from our experience to allow for an algorithmic symbolic system that, in turn, allows one to more easily develop proofs, i.e. deduce logical consequences of the principles, with much-reduced chances of error. By logicizing and symbolizing, one creates an algorithmic, procedure-orientated approach that also facilitates discovery. This is the main reason modern mathematics is so much more powerful and advances so much more quickly than did ancient mathematics.[9] It is also the reason math's meaning is often misunderstood.

# IV. Axiomatization and Extension
# of the *Fundamental* Numbers

To begin to understand the interesting depth that is floating all around modern math, let us see how this axiomatization can be done with the principles we discussed.

---

[8] Indeed, as we will see, in modern mathematics, "one" is used in two analogical ways: a simple unit and a complex unit (a fundamental number). Two is a *whole* number; concentrating on its holicity, we can say two is "one" in a different (analogical) sense, such as when we say 2/2=1. In this equation, we measure two by itself. Of course, the *plurality* of a thing (its nature as something with multiple parts) is distinct, though not separate, from its *unitive* aspect (the fact that it is a whole). This is the essence of the real difference (not merely definitional difference) between *number* and *one*.

[9] Another reason is physics constantly stretches us to see and incorporate wider principles into our mathematics.

Some principles, most notably the concept of number, are primary (because of their generality) and thus will continually be extended by analogy to be used in most all that will come later. The names we use in high school and college mathematics reveal the primacy of the analogical concept of number: natural numbers, rational numbers, irrational numbers, real numbers, imaginary numbers, complex numbers, quaternions which are aptly called by one teacher "hyper-complex" numbers,[10] etc. Getting this beginning right is essential to getting everything else right.

Recall the essential aspect of number is that *it is a whole whose parts do not share a common boundary and is completed by a last part*. Jumping off the latter point, after any given whole with parts, i.e. number, there is another whole with parts completed by a new unit. This is the concept of successor, which we will incorporate into our axioms once we have "broken up" the principles of number into atoms, i.e. small intuitive pieces. Another piece or aspect of number is the whole that has no parts, i.e. the simple unit, which we designate as "one" or "1".

Any given fundamental number has parts. Keep in mind that, because of this, the unit (i.e. one), that *by which* we "measure" the whole and which has no parts, is not a fundamental number. Recall that the meaning of unit is pulled from the real; it is not something we simply invent -- though we will axiomatize and create a system for this bit of reality along with others. We will call this aspect that we have isolated the "principle of the unit." Like the other principles, we will need to incorporate it into our axioms.

Note, however, that the propositions and symbolic structures that we will enumerate and set up will not have to convey *all* the aspects of the original concept. Indeed, they should not, for this would be an indication that the very process of symbolizing of the real and reflecting it into logical relations (logicizing) that we seek to accomplish has not been accomplished. In short, we would not have ignored enough of the complex distinctions within the thicket of the real to make a basic system of symbols subject only to basic logic. We need to "flatten" the real to make proving and calculating as "automatic" and "mechanical" as possible so as to preserve our ability to grind smoothly ahead with no error--or certainly much less chance of error. Later, we can bring back some of these distinctions in different axioms (at the possible expense of others) as we choose to explore different aspects of the real quantity.

We now have all the basic components needed and thus are ready to introduce the axioms of the *natural numbers*[11], *symbolized by N (no F subscript)*; these are the first type of numbers introduced in modern mathematics texts. Giuseppe Peano (1858-1932) axiomatized the natural numbers in the following way (I comment on the incorporation of the primary concept of number after each axiom.):

1. One is in *N*.

> *We flatten the distinction between that by which we measure (which results in 1 or one) and that which we measure (which results in $N_F$ or fundamental numbers) so as to treat all elements of N with the same symbols and manipulation rules.* **One is now called a number** *by ignoring its distinctive character to some degree. Obviously, we cannot get off the*

---

[10] http://www.lboro.ac.uk/departments/ma/gallery/quat/intro.html

[11] Note our definition of "natural number" $\left( N \in \{1, 2, 3...\} \right)$ is common in number theory, but in set theory zero is usually included as well.

      *ground without that by which we measure in some operational sense; this axiom begins to bring it in.[12]*

2. All $n \in N$ have a unique successor $n'$ in $N$. (Note: the prime symbol means successor), so that, for instance, $2' = 3$ *and* $3' = 4$.)

      *Numbers in the primary sense (i.e. not the unit, "1") are completed by a last unit, also called a successor, thus a new number can always be had.*

3. a. There is no $n \in N$ which has 1 as its successor.

      *Principle of unit: i.e., incorporates some of the uniqueness of the unit.*

  b. Each other $n \in N$ (i.e., $n \in N_F$) is a successor of some $m \in N$ (i.e., $m' = n$).

      *Since number is completed by last unit, there is no other way to get a number. Without this, our system would not "know" that there are no other numbers than these.*

4. Every distinct natural number has a distinct successor,
i.e., $m, n \in N$, $m \neq n \Rightarrow m' \neq n'$.[13]

      *Completing the previous units by a last unit makes each number unique.*

5. Principle of induction: if a statement is true for 1 *and*, given it is true for $n$ and is also true for $n'$, then it is true for all $n \in N$. Said another way:
      If $C$ is a collection of natural numbers (i.e., $C \subseteq N$) *and* 1 is in that collection ($1 \in C$) *and* $n \in C$ implies $n' \in C$, *then* the subset $C$ is $N$,
      i.e., we automatically get all of *N*.

      *Again, in seeing the concept of number, we see that it entails that a new number can be had by simply adding a unit to a number that I already have. Induction, this last Peano axiom (# 5), captures in a systematic way the potency of number to continue without limit by successive ones. Indeed, the principle of induction captures the essence of number.[14] For, a statement that meets the requirements that it works for one and its successor means that it meets the essential requirement of number and thus must be true for all numbers.*

The beauty of abstraction (leaving out specifications) is that we can say things that are true in many different incarnations in the physical world *and* that we can say general things about whole classes of objects and classes within classes, etc.

    The process of axiomatizing the natural numbers (making a symbolic logical system that incorporates the principles) quickly leads to the concept of the algebraic field. To lay bare the connection with the real that is our task in this article, we will not build on these axioms, but will develop the axioms for a field from direct sensorial experience. In fact, the Peano axioms are specific to the natural numbers and thus include too much specific information to allow defining the general concept of a field or even an ordered

---

[12] A crucial point here is that because "one" is intimately related to the numbers (one is that by which they are measured) we can make a system of axioms in which that relation is embodied and talk about these entities as if they were one thing by focusing on their interdependence.

[13] Or equivalently, the contra-positive is: $m' = n' \Rightarrow m = n$

[14] It does so by implicitly assuming the principles incorporated by the axioms before it.

field. We thus start from what we know of numbers, taking advantage of the insights already accrued.[15]

# V. The Roots of the Axioms of High School Math

## A. Meaning of Addition and Multiplication

One of the first things we notice about numbers is that any given fundamental number can be decomposed. Obviously, we know already that each number can be decomposed into units. But, less trivially, larger numbers can be considered parts of another even larger number. This is fundamentally the concept of *addition*. In the real world, before we have abstracted everything but number (parts outside of each other that do not share a boundary), we can think of addition as the action of adding new parts not already there. Still, the mathematical entities are not changing; they are what they are; 2 is 2, and not 3 and does not ever cease to be 2. They are static. Hence, for instance, $2 + 2 = 4$ is most properly thought of as a relation between 4 and 2. Of course, it is always helpful to have the active presentation (as if we are changing 2 into 4) somewhere in the background ready in our imaginations to aid in thinking about addition and the other "operations." After all, as children, before we knew what 4 was, we played around with things, say marbles in a bowl, until we saw four has, in a way, two twos in it.

A minimal statement of addition, i.e. adding just **two** natural numbers ($a + b = c$), manifests the relation between three numbers. Of course, addition can also be thought of as having as many arguments as you please ($a + b + c + d + \ldots = x$). In particular, consider the statement $1 + 2 + 3 = 6$; this is revelatory of the nature of 6, 1, 2 and 3, in that they can be combined in 6.[16]

Of course, we have, in our understanding of addition, introduced a whole with parts, which are, in turn, also wholes with parts. Indeed, we have implicitly introduced fractions. With this in mind, we might rephrase our above statement as: $1/6 + 2/6 + 3/6 = 1$. Or in words: we have one part (which is a simple undivided part, i.e., a unit) of a whole 6 and another part (which has 2 parts (units) to it) and a third part (which has 3 parts (units) to it). This sum equals one, but it is not the simple unit, not a simple "1", but a composite "1"; indeed it is 6, and at that a special 6; its parts are not all units but have parts themselves.

Again, in our current axiomatization, we want to ignore this aspect of parts with parts and thus import key aspects of quantitative reality into our system without *having* to make those detailed distinctions that would sacrifice the simplicity and abstractness that we can get at this level.[17] So, *addition* is really a way of seeing the composition of a number, a whole with parts. Or recalling the background active picture in our imagination (e.g. adding marbles to the bowl picture --increasing parts of the whole)[18] from which we first got the concept of addition, we can say addition tells us what numbers go together to make another number.

---

[15] There is, of course, much more to be said of the fundamental concept of number. These axioms do not lay out the specifics of each number, but only what is generically true of them all. Along with the generic truths about any given number, each number brings something new; for instance until we analyze 49, we do not know it is composed of seven sevens. Number theory studies the many fascinating aspects of number, including the fascinating *Fundamental Theorem of Arithmetic* that every number can be composed from a unique multiplication of primes.

[16] Notice that, here, in order to define addition in a useful way for our system, we treat one as a number (i.e., as one of the entities represented by a symbol to be discussed and manipulated) and draw no direct attention to its unique character.

[17] Though later, we may want to include it and move to another level.

Now, *multiplication* is a special type of addition. Most properly, i.e. in the most primary sense, it applies to those numbers that allow multiple numbers "within" them. For instance, four has two twos. Eight has two fours or four twos ($8 = 4 + 4 = 2 + 2 + 2 + 2$). So, in multiplication, we have a number added multiple times to itself. We have a number (whole with parts) that is composed of more than one of the same number.

By analogy, we also regard "taking one of a number" as multiplication; i.e. $1 \times 4 = 4$. But why do this? After all, the statement is fundamentally the same as simply considering the number itself. We do so to complete the logic of our system. Similar to the way we sometimes flatten reality so that our system concentrates attention on a given aspect (thus not having to consider the complication of other aspects), we, here, stretch a reality to establish the system. Our system simply takes advantage of the obvious fact that a thing is related to itself as itself and thus can be called a relation in a logical sense and thus, by analogy, lumped with relations between different numbers.

In $4 \times 1 = 4$, we are noting that there are 4 units in four, but multiplication in our system is defined for "numbers," so we here make implicit use of our generalization of "number." In this way, we incorporate in our definition of multiplication 1) the fact that the *simple unit* is that by which we measure fundamental numbers *along with* 2) the fact that some larger numbers can be measured by multiples of certain smaller *fundamental numbers*, but we do so *without making the distinction* between these two facts *explicit*.

*Subtraction* is just a different perspective on addition. Instead of noting that $1 + 4 = 5$, we note $5 - 1 = 4$; both express the same relation between the three numbers but each emphasizes a different aspect. In the former, we note that a last unit completes the 5 starting from four. We note the same in the latter. The difference lies in the *background* active viewpoint. The latter emphasizes how, if we go back to an image in the mind, leaving the abstract notion of number for a moment, we can make four from five; say by taking one marble out of the center pile of five.

Similarly, *division* is another perspective on the special type of addition that is multiplication.

Beyond understanding addition, subtraction, multiplication and division, we need a symbolic framework that is more compact and precise than words to carry the meaning and logic of our thinking about number.

## B. Symbol Usage

We define the letters "*a*" through "*d*" to be symbols that represent one of the natural numbers (1, 2, 3, 4…), i.e. $a \in N$. Recall we include "1" for the reasons we have already explained; it is an analogical use of the word number that we point to by calling them natural rather than primary or fundamental. In our manipulations and formalism such distinctions are left out, but that does not mean they are not there and that we should not be mindful of the full reality at some point.

Using the symbols "*a*" to "*d*," we next define the rules of manipulation called the reflexive property, the symmetric property and the transitive property. These rules capture some of the properties of the naming involved in symbolization. Symbols always involve naming of a sort. The name "Fred," for instance, is a symbol by which we refer (say) to the man who lives next door. Similarly, the definition of "=" and the meaning we give to the variables "*a*" and "*b*," i.e. as purely names of some unspecified mathematical object, result in necessary properties of these objects as names.

---

[18] Or more profoundly in the embryo subdividing picture, the embryo manifests more clearly the idea that a number is a *whole* with parts.

In particular, names, in order to be names, must have the following properties. First, the *reflexive property*, written as "$a = a$," simply recalls the self-evident fact that every name must always be itself; obviously, everything must be itself. Of course, we could endow the placement of the name on the right or left of the equal sign with a meaning, but we choose to make the designation of the name come only from the letters (and their order in a word) of which the name is composed and not sensitive to where they are on the page. This same lack of use of left-right position in naming is involved in the *symmetry property*, written as: $a = b \Rightarrow b = a$. Symmetry encapsulates the simple fact that when two names are assigned to the same object, they are assigned to the same object in both contexts. The last property, the *transitive property*, written $a = b \text{ and } b = c, \text{ then } a = c$, is the equally obvious statement that if two names are said to refer to the same object as a third names does, then all three names refer to the same object. [19]

We now have all the ingredients needed to prove certain properties necessary for manipulating natural numbers.

## C. Understanding the Laws of Addition

The essential standard axioms of addition are (with a note indicating where they come from in reality):

<u>*Commutative Law of Addition*</u>: $a + b = b + a$

This is obvious if we recall that the parts of any number are not in any particular order (they do not share any common boundaries to distinguish them); hence, the relation to the numbers or units that compose the given number is not dependent on the order of consideration.

<u>*Associative Law of Addition*</u>: $(a + b) + c = a + (b + c)$.

The parentheses[20] here simply mark the order in which we consider the various numbers, which, as we have just pointed out above, makes no difference.[21]

---

[19] We can also add another axiom to make explicit the substitution principle in equations, which is simply again a feature of multiple names for the same object. The axiom, called the *Substitution Axiom of Equality* can be stated as: $(a = c \text{ and } b = d) \text{ implies } (a + b = c + d \text{ and } a \times d = c \times d)$.

[20] In introducing parentheses, we consider addition as primarily a binary operation, having only two operands; this reduces the complex consideration of *all* the relations of the various numbers in the sum to each other to considering only two at a time. Those two are then seen to be part of one number and reduced to that number in order to compare that number to another number. These latter two numbers are, in turn, seen to be part of a third number (in the sense previously described). So, for instance, (3+2) + 4 , becomes 5 + 4 = 9.

[21] If we were to "leave in" the fact that numbers composed of other numbers are qualitatively different from the most primary concept of number, we might make a distinction, but we have chosen not to. When we define rational numbers, we will reintroduce this fact through the back door.

The essential standard axioms of multiplication can now be shown to be true in the following way. By "true," we mean in accord with the reality of number that we originally "pulled-out," or abstracted from physical reality.

## D. Proof of the Laws of Multiplication

<u>*Commutative Law of Multiplication:*</u> $a \times b = b \times a$

| | | Reason |
|---|---|---|
| $a \times b$ | $= b \times a$ | |
| $\underbrace{b + b + b + ...}_{a \text{ times}}$ | $=$ | *definition of multiplication* |
| $a + \underbrace{(b-1) + (b-1) + (b-1) + ...}_{a \text{ times}}$ | $=$ | *definition of subtraction* |
| $a + a + \underbrace{(b-2) + (b-2) + (b-2) + ...}_{a \text{ times}} =$ | | *definition of subtraction* |
| $\underbrace{a + a + a + ...}_{b \text{ times}}$ | $=$ | *definition of subtraction* |
| $b \times a$ | $= b \times a$ | *definition of multiplication* |

*Q.E.D.*

The above can be illustrated using the following picture sequence shown in figure 1a. Recall, however, that the *concept is bigger, much more general, than the picture sequence used to illustrate it.*

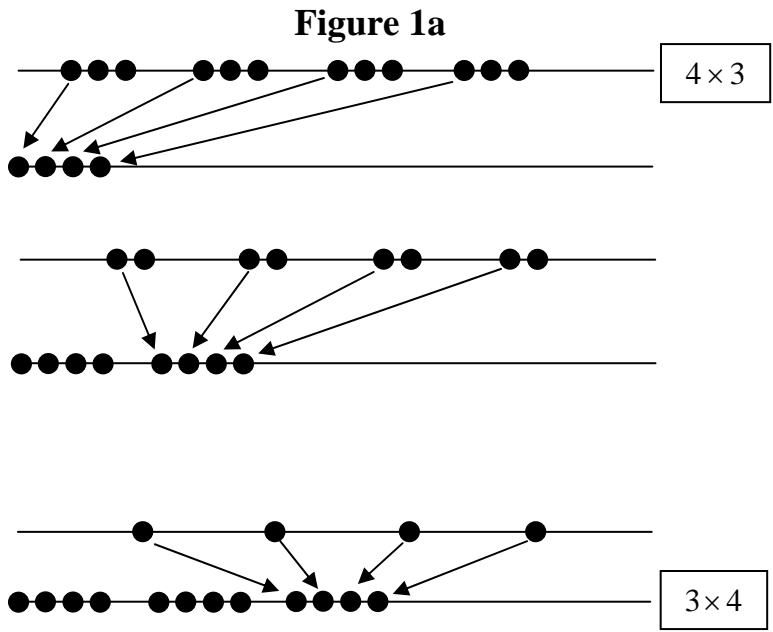

**Figure 1a**

$4 \times 3$

$3 \times 4$

## Figure 1b

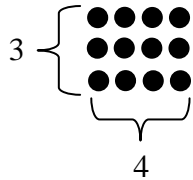

It is hard to overemphasize that "number" is a reality in things, not existing separate from them but in them, and we abstract it *from them*. In the abstract concept of number, we have something that covers more than any particular realization we may use as an example. We thus cannot use the diagram as a proof, for the diagram strictly as the particular diagram is true only for the particular number shown in the particular circumstance shown (e.g., for Figure 1b for two dimensionally extended bodies).

The abstraction allows us to talk about all numbers in all circumstances. This dual (abstract and real) aspect of mathematics is critical. Math is both general, whence it is awesome power, *and* real (having come from the real), whence its awesome ability to shed light in physics.

The proof for the associative law is:

*Associative Law of Multiplication*: $(a \times b) \times c = a \times (b \times c)$

| | | |
|---|---|---|
| $(a \times b) \times c$ | $= a \times (b \times c)$ | |
| $\underbrace{(b+b+b+...)}_{a\ times} \times c$ | = | *definition of multiplication* |
| $c \times \underbrace{(b+b+b+...)}_{a\ times}$ | = | *commutative law of mult.* |
| $\underbrace{\underbrace{(b+b+b+...)}_{a\ times}+\underbrace{(b+b+b+...)}_{a\ times}+\underbrace{(b+b+b+...)}_{a\ times}+...}_{c\ times}$ = | | *definition of multiplication* |
| $\underbrace{\underbrace{(b+b+b+...)}_{c\ times}+\underbrace{(b+b+b+...)}_{c\ times}+\underbrace{(b+b+b+...)}_{c\ times}}_{a\ times}$ | = | *commutative prop. of addition* |
| $\underbrace{(c \times b)+(c \times b)+(c \times b)+...}_{a\ times}$ | = | *definition of multiplication* |
| $a \times (c \times b)$ | = | *definition of multiplication* |
| *So*: | | |
| $a \times (b \times c)$ | $= a \times (b \times c)$ | *which is commutative law of multiplication* |
| *Q.E.D.* | | |

Finally, we have the combined axiom of the distributive law:

<u>*Distributive law of multiplication over addition:*</u>   $a \times (b+c) = a \times b + a \times c$

$$
\begin{aligned}
a \times (b+c) &= a \times b + a \times c & \\
\underbrace{(b+c)+(b+c)+(b+c)...}_{a \text{ times}} &= & \text{\textit{definition multiplication}} \\
(\underbrace{b+b+b+...}_{a \text{ times}}) + (\underbrace{c+c+c+...}_{a \text{ times}}) &= & \text{\textit{commutative prop. of addition}} \\
(a \times b) + (a \times c) &= a \times b + a \times c & \text{Q.E.D.}
\end{aligned}
$$

We can, of course, also prove the above laws by starting with the Peano axioms. The gain in so doing would be more ability to automatically manipulate symbols and thus less possibility of error in deductions in long proofs. On the other hand, in that approach our connection with the real is thinner and more obscure-- logic and symbols being in the foreground instead of the meaning to which they refer. Indeed, the reduction to steps in logic executed through manipulation of symbols is the very reason the manipulation is more mechanical (and thus readily put on a computer).

By contrast, the above proofs are much more easily accessible to the non-mathematician. They are based on our immediate experience and can be learned by very young children. My 7 year old son was able to master these proofs, and give, on his own, a detailed explanation of the first to his grandparents. There is much to be gained by emphasizing the principles and where they come from before we introduce the axiomatic systems that capture their deep, complex and even beautiful meaning.

Notice that these proofs are, strictly speaking, only valid for the natural numbers, {1, 2, 3…}, though they can be fairly easily extended to include the more general meanings of number that we discuss below.

One more statement can be added at this level, which can also be readily generalized. It is the *axiom of closure* of multiplication and addition; i.e., when we add natural (or fundamental) numbers, we always get another natural number. Why so? The reason is clear when we recall that addition is nothing but the relation of three numbers: the two summands and the sum. A little thought about the nature of a number shows this is closely related to the principle of induction.

That principle, as we saw in Peano's postulates, says that if a statement is true for: a) one and b) if assuming it is true for *n* then it is true *n* + 1, then it is true for all numbers. This principle, as we also saw, contains the essential aspects of number, including that every fundamental number is terminated by a last unit. Hence, the closure of addition follows from the fact that, given that any two numbers are made up of successive units, there must be a third number corresponding to all the successive units taken together terminating the whole by a last (arbitrary) one of all those units of the two numbers together. Again, this is clear from the concept of number we get from the physically real; consider, for instance, 52 cards from a playing deck (the whole) given evenly to two people, or given 50 to one, two to another.

## E. The Axioms for Integers

Now, the above properties and laws (namely, the reflexive, symmetry and transitive properties, the commutative and associative laws of addition and multiplication and the distributive law) apply some of the properties of the fundamental numbers. They do not yet incorporate important properties of the natural numbers. We need, for example, to add the properties of one, which were implicit in the Peano axioms but not yet in the above system. After incorporating some aspects of one, we will add zero and introduce the negative numbers to make a system that applies to all the integers.

To include the properties of one in our system of axioms, we require that for all $a$, $a \times 1 = a$. For this reason, we call "one" *the identity element for multiplication*. The axiom of the multiplicative identity element states that there is a unique multiplicative identity element.

Note again we use a symbol, e.g. "$a$", to refer to the mathematical entity that we are trying to specify properties for in our axioms. Further, this entity to which "$a$" refers is called a number in a (analogically) general sense. We sometimes specify what type of number by adding an adjective such as *natural* number or *rational* number or *real* number. It is sometimes said that "number" itself is undefined in our axioms, because we are concentrating on symbolizing and logicizing the essential properties, not understanding them *per se*. Of course, ultimately, the entire reason we do this is to understand them; it is just that we usually forget this amidst all the work and effort in proving and using the formalism, which still keeps us in the presence of the properties *indirectly* in the powerful system of axioms and theorems.

To include zero, we require, for all possible "$a$": $a + 0 = a$. Zero is called the additive identity element. The *axiom of the identity element for addition* states that there is a unique additive identity element, $0 \neq 1$. Here we are simply (though you may say somewhat oddly) extending our definition of addition to include the relation of any number, a, to nothing. Namely, in a way that reveals that this axiom incorporates more logic than mathematics, we can say "$a$" is composed of "$a$" units and also nothing. Of course, something is not really composed of nothing, but we can by analogy include it in our system. We can also think of it in the "active" constructive sense of doing nothing to "$a$" units. The idea of zero, which represents an absence of anything, also prepares the system of axioms for what comes next: negative numbers.

Indeed, in the active way of thinking about subtraction, we can say if I have "$a$" things and someone takes all "$a$" of them away, I have none left; in other words, I have the absence of those "$a$" things that used to be there, i.e., nothing. Symbolically, $a - a = 0$. Now, look at this equation more accurately by thinking about the subtraction as the relationship among three numbers; i.e. as we pointed out before, look at subtraction as a different perspective on addition. In so doing, re-write the equation as: $a + (-a) = 0$. We say "$-a$" is simply a new, extended yet again, definition of "number." We then can say that "$a$" and "$-a$" are simply the two "numbers" that go to make up the "number" zero. This gives the generality we need, for it works for all integers, whether positive, negative or zero. Succinctly, this latest axiom, the axiom of the existence of additive inverses, states that for all numbers "$a$" there is a number "$-a$" such that $a + (-a) = 0$.

This axiom can be thought of as very directly incorporating the idea of "taking away" in our formalism, but more importantly it closes the loop on subtraction, rounding out our system of axioms with respect to addition. For instance, before the axiom of additive inverses, what would we do with, for instance, $6-8$? With the axiom, we simply write $6+(-(2+6))=-2$.

Furthermore, inclusion of the additive inverse axiom does not affect our closure axiom for addition, for, as we have said, we have simply included the same numbers with a different background interpretation, namely that of taking away, which means we have really added nothing new except to make more robust the various relations that numbers have to each other in their composition. Indeed, the very definition of the inverses makes manifest that no new elements can be made by defining negative numbers. More importantly, if we had not already included zero, we would have to now in order to complete the logic of our system. Of course, the previous inclusion of zero through the identity element of addition cannot affect closure in any other way, for by that very axiom it does not alter the relationships among numbers; in this sense, we could certainly say zero adds nothing.

Now, what about closure with respect to multiplication? To answer, we have to say what we mean by multiplication of negative numbers. For instance, what does $-2\times 2=-4$ mean? First, note we are, as above, thinking of "$-2$" as an entity, as if it could exist as such. As such, "$-2$" is a *purely mental construct* necessary for our axiomatization, but it is reductively real; its parts are real. We build the idea of "$-2$," from the real fact that two parts of a whole can be thought of as missing or present. In the background active point of view, we can imagine an object that used to have two parts or from which two parts are taken away, leaving none. The minus connotes the background operation that is simply a way of incorporating the fact we can have two of something or not. So, $-2\times 2$ can be read as taking away two twice, or in other words as taking away the number that it implicitly refers to, i.e., 4; after all, we know the relation of 2 to 4 is 2 + 2 = 4. In other words, $-2-2=-4$.

How about the case: $-1\times(-2)=2$. Here we are taking back *or undoing the fact that we took 2 away*. This is simply a convoluted way of saying we have 2. Consider then $-2\times(-2)=4$. This kind of multiplication, as we have been implying, is a special case of subtraction: $2\times(-2)$ means subtract 2 twice: i.e., $-2-2$. While, $-2\times(-2)$ means subtract $-2$ twice or $-(-(2\times 2))$, which, in turn, means: take the number four away, then take away that we took it away, yielding 4. In this way, we incorporate simple facts about the potentiality of things to be there and not be there with their nature as numerable into our system.

With this, we see clearly that we do maintain *closure with respect to multiplication* with our extended definition of number that includes zero, one and the negative numbers. Namely, for any two numbers, "*a*" and "*b*," defined as above, multiplication, $a\times b$, will always yield a number of the same type.[22]

*The axioms that we have so far enunciated plus the axiom of the multiplicative inverse (which we will discuss in the next section) are the system which defines a **field**. Any mathematical entities that satisfy those axioms will be called a field. **Ordered fields***

---

[22] Viewing the numbers as elements of a group, with the binary operation of multiplication, we would express the closure of the group by saying: for every entity *a* and *b* in the given group, $a\times b$ is also in that group.

import more information about the fundamental concept of number into our system. To have an ordered field, one adds the following axiom and properties[3]:[23]

*Order Axiom of Comparison*:
Exactly one of the following is true: $a > b \text{ or } a = b \text{ or } a < b$
*Order Transitivity Property:* $a < b \text{ and } b < c \text{ implies } a < c$
*Order Property of Addition:* $a < b \text{ then } a + c < b + c$
*Order Property of Multiplication*: $a < b \text{ then, for } c > 0 \quad a \times c < b \times c$
*Substitution Property of Order*: $(a < b \text{ and } a = c \text{ and } b = d) \text{ implies } c < d$

Since these properties follow from a similar type of analysis that we have done above multiple times, we simply list these and leave the analysis showing the connection with the real to the reader.

With this last set of axioms (excluding the axiom of the multiplicative inverse) we have a system that captures many of the properties of integers. We denote the integers by the symbol **Z**, $Z \in \{Z^-, 0, Z^+\} = \{...-3, -2, -1, 0, 1, 2, 3, 4...\}$; note that $Z^+ = N$ and $Z^- = -N$. This is, then, our third analogical use of the word "number": first we had the fundamental numbers, then the natural numbers and now the integers. Considering how often we have used analogies, it is starting to look like we mathematicians are more like poets than we might have thought. Unlike poetry, which uses metaphors that are loose types of analogies, the analogies in mathematics are very tightly connected, and we are here precisely delineating the meaning. Also, unlike poets, we, in our axiomatic systems, flatten our analogies and do not keep them up front, but quickly leave behind the distinctions that separate the various types of numbers, so as to bring to the fore what is common. This allows us to move swiftly forward in discovery and proof, without being so prone to erroneous leaps in logic, though we will sometimes, even often, leap over our understanding of how we get where we end up and where it all came from.

## F. Axioms for the Rational Numbers: the First *Fields*

As in our previous generalizations of "number," we include more information from the real to expand our system to the more general type of entities called the rational numbers. We go back to our insight that addition is a relation between three numbers, in which we think of the possibility of numbers inside of other numbers. We can have parts that are not simple, but which themselves have parts. A rational number is a ratio of integers or a relation between integers (e.g. 1/2 or 16/95).

From an axiomatic point of view, we also note that it would be nice to completely "close the loop" with respect to multiplication by making division a closed operation, in analogy to the way that we made subtraction a closed operation by adding negative numbers. We, thus, say that, for all "*a*" and $b \neq 0$ in the mathematical entities ("numbers") that we are now defining, *a/b* is also one of those numbers. Or said another

---

[23] Another way to capture the key aspects of these axioms, which more closely links them with the natural numbers, is to use instead the following axioms of order: For elements $x \in$ {entire set that satisfies the system of axioms}, which is the real numbers, and for elements $a, b \in S$, where *S* is a subset of the natural numbers (or more generally a subset of the positive reals) the following are true: I. $(a+b) \in S$  II. $(a \times b) \in S$  III. $-a \notin S$  IV. $x = 0 \text{ or } x \in S \text{ or } -x \in S$. (see [4])

way, for every $a \neq 0$ there exists an inverse "b" such that $a \times b = 1$. This is called *the axiom of multiplicative inverse*. "b" can be written as $1/a$.

With this statement, we have now expanded the entities described by our formalism to the rational numbers, which gives the smallest ***field*** containing the integers. Furthermore, before introducing this axiom, we had not defined division in our axioms.

We call them rational *numbers* in analogy to our primary concept of number, because they can be thought of as a whole, though in fact they are the relation between two numbers. For instance, 1/2 is the relation between the unit and two, the whole which has two parts (that share no common boundary). Similarly, 2/4 is the relation between 2 and 4; of course there is a proportion between 1/2 and 2/4. Indeed, in terms of *the relation between the relations*, they are the same; that is, 1 is to 2 as 2 is to 4 or 1/2 = 2/4 in the familiar way. Similar statements can be made for the negative numerators and denominators, as well as with numerators of zero.

More to the point of our axiomatization, any number in this extended sense can be thought of in this relational way. 2, for instance, can be thought of as 2/1. Still, 2 is different from, say, 2/4, for 2 is a whole with unit parts whereas 2's relation to four is as parts with units within it. But, similar to what we have done before, we ignore that difference to concentrate on the aspect that is the same. These flattened analogies are incorporated into our logical system of axioms and symbols, thereby allowing us to carry around and manipulate analogically general properties as though they were truly general. In this, we extend our definition of number, now being able to lump all previous numbers under the heading of rational numbers, so that fundamental numbers, natural numbers and integers in addition to those ratios not reducible to integers are all called rational numbers. Different though they are, their similarities are used to unify them in the axioms of rational numbers.

Having included all the ratios of numbers as numbers, i.e. entities that obey the axioms, we now may wonder whether multiplication is still a closed operation. Namely, will multiplication of any rational number by any other always yield another rational number? To answer we must first define what we mean by multiplication of rational numbers. We have said multiplication is a special type of addition, that is, numbers composed of multiples of the same numbers rather than the more general case which can have multiples or not. Consider: $2 \times \frac{1}{4}$. Clearly we mean by this: $\frac{1}{4} + \frac{1}{4}$. However, what if we have a fraction times something, such as: $\frac{1}{2} \times 4 = 2$? This is also clear. Though it no longer is a special type of addition in the most explicit way, it is in the deepest way. It gives the relation between 1, 2, 4 and 2. Why do we repeat the 2? Though the symbol "2" we use is the same, we want to distinguish one 2 from the other. A figure may help.

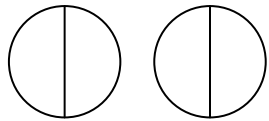 **Figure 2**

Note, first, that the circles in figure 2 have two parts; they are wholes with two simple parts. This is the first meaning of 2. Next, there are two circles. From this, we get a different type of 2; this 2 does not have simple parts but, instead, parts which also have parts.

In summary, we see in this 4 that there are two 2's (and each of those 2's has 2 ones). This 4 clearly includes the type of wholes with parts that are themselves wholes with parts that rational numbers always invoke. The same is true for $\frac{1}{2} \times \frac{2}{8} = \frac{1}{8}$ with respect to the 1, 2 and 8; it says that, considering two as a complex part within eight, we can always consider one of the two parts within the two.

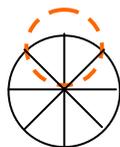 **Figure 3**

Now, $\frac{1}{4} \times \frac{1}{5} = \frac{1}{20}$, because "$\frac{1}{4} \times \frac{1}{5}$" means break each of the 5 parts of a given whole (i.e., "1") into 4 parts each. Or from the more proper perspective of static mathematics: consider a whole with twenty parts. Considered it as composed of five parts, each of which has four parts to it, so that any one part can be considered a unit in a full 20. Considering instead $\frac{3}{4} \times \frac{1}{5}$ simply means that we add this result three times, i.e. $3 \times \left( \frac{1}{4} \times \frac{1}{5} \right)$, giving 3/20. In this way, we arrive at the rule to multiply the numerators and denominators to get the new fraction. Of course, in some special cases we can, by using the equivalence of ratios, reduce the fractions in familiar fashion.

Hence, it is clear, by our definition of a rational number as a ratio of two integers (symbolically: $Q = \{all\, \frac{i}{j},\, where\, i, j \in Z\, and\, j \neq 0 )$ and by the definition of the multiplication of two rational numbers as simply the product of the appropriate numerators and denominators, that such multiplication will always yield another rational number. That is, closure with respect to multiplication is maintained. Again, note this would not be the case for integers with respect to division. Obviously, in dividing 3 by 4, we would not get an integer; in a way, this fact of the non-closure of division would drive us to generalize to the rationals. This is another example of how logicization and axiomatization can be very powerful tools leading one towards new mathematics--this despite the fact that it comes at the expense of intuition.[24]

Now, the new definition of number as a ratio of integers makes it clear that every such number will have a multiplicative inverse, which when multiplied by this inverse will yield the identity element, 1. Note, this "1" is our analogical *extension* of the unit; it no longer always means just the simple unit, i.e. a whole with no parts, but can also mean a whole with parts such as $\frac{4}{1} \times \frac{1}{4} = \frac{4}{4} = 1$; in words, there are 4 parts in 4. As an example, a chicken embryo that has divided four times is still one. Or, if I eat four pieces of a pie that has been cut into four pieces, I have eaten the whole pie.

---

[24] For those who need further exposure to how symbolic logicization can thwart intuition, there is a proof floating around on the internet which purports to use induction to prove that everyone is the same age. The site organizers do not actually believe this; yet, they challenge you to find the flaw in the train of reasoning. It is not a trivial task, even though no one would think for a moment that all people are the same age. (cf. http://www.math.toronto.edu/mathnet/falseProofs/sameAge.html)

## G. The Number Line:
## Moving Back to the Continuous (to Geometry)

Notice that, up to this point, each time we have extended the meaning of number, we have included more of the number line. That is, we can make a correspondence between fundamental numbers, natural numbers, integers and rational numbers and a line segment in the following way. We can assume that a line is divided into simple parts of equal size (this then being the unit by which we measure the line), and thus get a number to put at the end of each line in standard fashion. These are the natural numbers (or the fundamental numbers –number in the primal sense-- if we allow ourselves to notice, "1" is qualitatively different than the others).

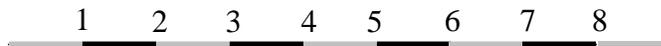 **Figure 4**

The integers can be represented by adding zero to the left on the above line and then extending the line to the left, symmetrically labeling that extended portion with the negatives of the natural numbers. The rational numbers can be represented by points intermediate between the appropriate numbers. For instance, 4/3 means a *whole with three parts* and a simple part of another *whole with three parts*. On the number line, we will take our old unit "1" and divide it into new units, so that our "1" now has three parts; we do the same for the old unit that separates one and two. In other words, four (new) units which are parts of a whole three would be placed to the right of "0." This puts the mark for 4/3 one third of the original unit to the right of 1 (see figure 5). The same is true for 5/3 etc., so that we could easily fill the number line beyond what we could write in the space available. Indeed, the number of rational numbers *between each integer* is infinite, and, in fact, it is of the same order as the integers themselves! By "of the same order," we mean that one can assign an integer to every rational number that one could think of and vice versa; they can be put, in some sense, into one-to-one correspondence.

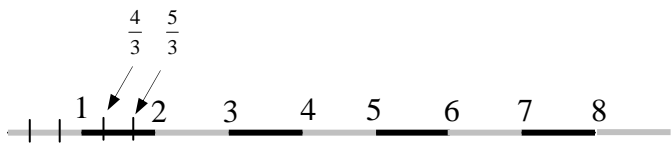

**Figure 5:** A number line illustrating some rational numbers. 4/3 is illustrated as 4 simple units, 3 of which make up what we call "1." That is, we have made new simple units so that "1" is no longer a simple unit but is now composed of three simple units.

Still, there are gaps in the line. How do we know? Not all points on the line are an integral number of "fractions of the original chosen unit" from the chosen origin. Said another way, there are "numbers," such as $\sqrt{2}$ and $\pi$, that cannot be written as the ratio between two integers. For instance, it can be shown that there is no commensurability between the numerator and denominator in the ratio of the circumference of a circle to its diameter ($\pi$). This is obviously not the case for rational numbers, wholes made up of units that are themselves part of a larger whole. Again, not only is there no clear necessity that all quantities be commensurate like this, we know cases, such as π, in which it is not; we call them irrational, or not rational. By this, we mean they are not commensurate ratios, not that they do not make any sense. Irrational numbers, such as π, ultimately come from geometry.

To see how irrational numbers come from geometry (and finally from our senses) consider the Pythagorean Theorem as proved implicitly (remember figures *per se* are concrete and particular, whereas a proof must be general) in Figures 6a and 6b below through the use of geometric figures.

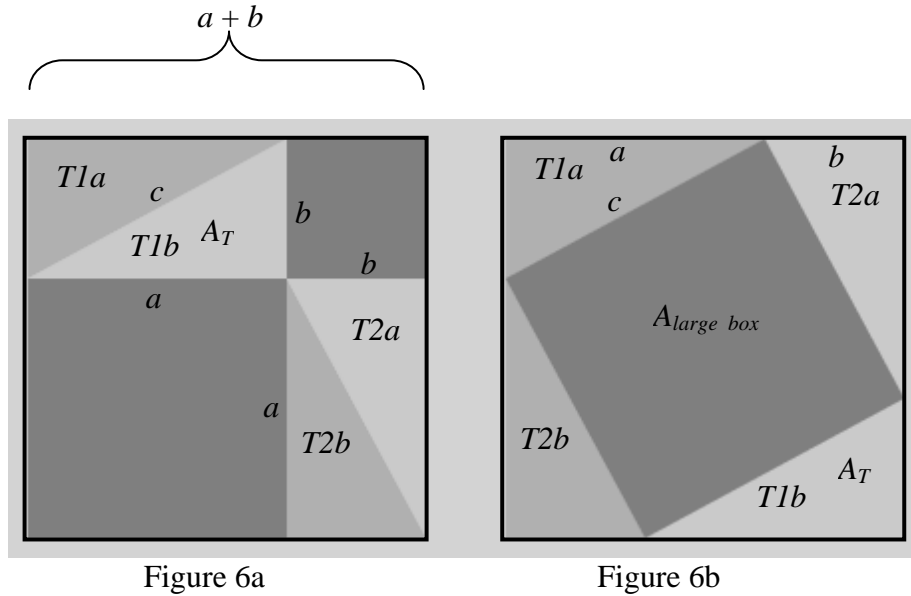

Figure 6a　　　　　　　Figure 6b

**Figure 6:** *To make Figure 6b from Figure 6a:* (1) erase the square with sides of length $a$ and also for those of length $b$. (2) Slide triangle *T2a* into the top right-hand corner preserving its orientation so that its short edge is pushed against the top side of the figure box. (3) Slide triangle *T2b* (keeping its orientation) over to left edge. (4) Slide triangle *T1b* into the bottom right hand corner preserving its orientation so that its short edge is pushed against the right side of the figure box. By careful construction, we can see that $A_{large}$ is a square. We also know $A_{large}$ is a square because its sides are equal, and the sum of the acute angles in each one of the four congruent triangles must be 90°, and thus all the angles of the $A_{large}$ are right angles. *Note that the area of overall box on left = $a^2 + b^2$ + 4 $A_T$ = Area of overall box on right = $c^2$ + 4 $A_T$. Thus, $c^2 = a^2 + b^2$.*

Notice that if we take $a = b = 1$ in the Pythagorean Theorem, $c^2 = a^2 + b^2$, we must have $c = \sqrt{2}$, i.e. a number such that if we add it to itself, itself times, then we get two. This we cannot do with rational numbers. We can prove this in the following simple way. Consider that there is such a ratio of integers, $i$ and $j$, in lowest terms such that:

$\sqrt{2} = \dfrac{i}{j} \Rightarrow 2j^2 = i^2$, thus $i^2$ must be an even integer, and thus so must $i$,

so take $i = 2m$ which gives:

$2j^2 = 4m^2 \Rightarrow j^2 = 2m^2$, which means $j^2$ must be even and thus so must $j$.

This means $\frac{j}{i}$ is not a fraction in lowest terms and this is a contradiction.

Hence, $\sqrt{2}$ is not a rational number. Therefore, if we want the full abstractive powers of number to apply to geometry, we must admit a new axiom to our "numbers," (thus expanding the meaning of the word). Or said another way, we introduce more of the structure of reality by a new axiom. In this case, we introduce some aspect of the continuity of geometry in the back door. This is how the intermediate structure between discrete quantity (fundamental number, also called arithmetic in ancient times) and continuous quantity (geometry) called *analysis* is built.

We do this by the axiom of completeness:
> ### *Axiom of Completeness*
> Every bounded, non-decreasing sequence of real numbers converges and its limit is a real number.

As an example of a non-decreasing bounded series, take the geometric series based on ½: $S_n = \frac{1}{2} + \frac{1}{4} + \frac{1}{8} + ..... + \frac{1}{2^n}$. We know it converges. A simple chug on the calculator shows that as we add more terms, i.e. as *n* gets bigger, we get closer and closer to 1. That is, the series $S_n$ can be explicitly written as: $\{.5, .75, .875, .9375, .96875...\}$.

The axiom of completeness serves to fill in the "gaps" in the rest of the number line. More precisely, before this axiom, there were infinitely many points (intersections with other lines) on our number line which were not associated with a number.

With this axiom, every point can be assigned a number. The numbers that obey the entire list of axioms (not including the Peano axioms) are called the ***real numbers.*** Obviously, they are an example of an ordered field.

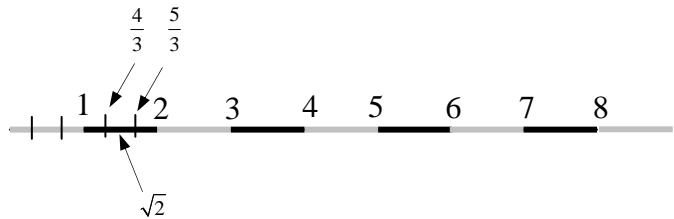

**Figure 7**

To see more closely how our final axiom fills in the gaps, consider the point associated with $\sqrt{2}$. We can approach the square root of a number by continually adding a sequence of rational numbers; there are algorithms that help you do this in base 10 that you may know. We then can get as close as we want to any given part of the line. Why? Again, because, though the numerator and denominator we need are not commensurate, we can get as close as we please by altering the numerator and denominator. We could effectively do this by adding smaller and smaller rational numbers as corrections to the given rational number approximation, [25] as we did in the geometric series.

---

[25] This can be best seen by the never terminating or repeating decimal representation of an irrational number. For instance, $\sqrt{2} = 1.41421356237...$. Other bases could, of course, also be used.

Indeed, we can, by a mental construct, act as if we have completed this addition and reached the required point and still call this limit, by analogy, a number. It is a number only in the extended sense that it obeys the axioms of number we have defined above. Again, generally,[26] such limits exist because of the continuity of the line, which we abstract from the real and import into our system via the *axiom of completion*. The new number obeys the rest of the axioms, because it is always arbitrarily close to a rational number.[27] Of course, the rational numbers are, in turn, related to our fundamental concept of number in the way previously described.

Hence, we have a system of axioms that gives us access to the continuity of the real line and, by an easy extension, to the continuity of any number of dimensions. Because the analogical concept of number lumps many qualitatively different types of "number" into one, not specifying the particulars or including anything that would bring out the differences, the axiomatic system that embodies such an analogical concept of number is a very powerful tool. As always, we should recall that the system of logic and symbols that we have put together gets its controlling features from real things from which we originally extracted its axioms. So, *ultimately,* any given system of mathematical axioms *is* about something, contrary to Russell's statement that we do not know what math is about or whether it is true.

# VI. Conclusion

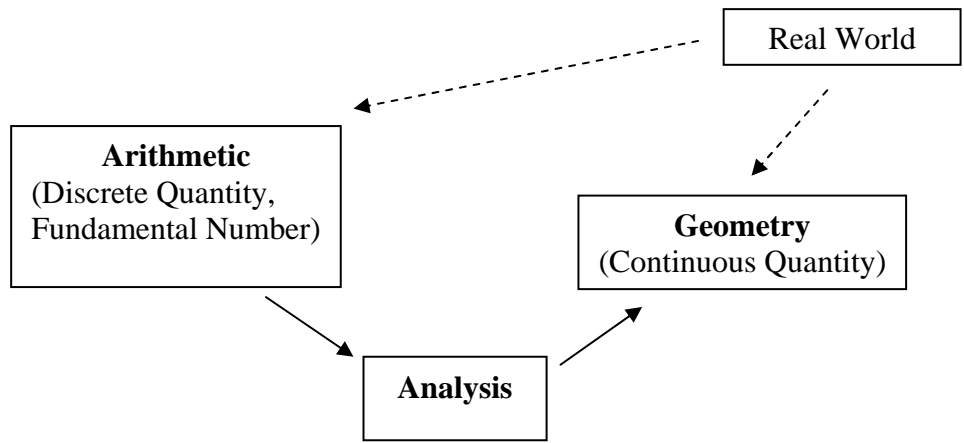

**Figure 8:** Chart of relationship between various mathematical disciplines

So then, what is math? Math is the study of quantity,[28] which we abstract from physical things that we see, hear, touch, taste and smell. Quantity is what we get when we leave all behind but a thing's extension. This leaves us with extended figures such as circles and squares and the like. From the parts of such figures, we can abstract, further, 23the idea of number. Hence, the study of quantity is divided into the study of continuous quantity, (parts outside of each other that share a common boundary, i.e. geometry) and the study of discrete quantity (parts outside of each other that do not share a common boundary). This latter subject, i.e. number, is the more abstract. Now, we can take the truths that we see in discrete quantity and import them into a system of axioms that will allow us to manipulate and carry around these truths. They also serve both as spring boards to understanding concepts implicit in them and as ground work for importing more of reality into our system. Often such axiomatic systems can appear to be the arbitrary creations of you or me, rather than facts about the world. This is only because, as educators and students of math, we have long habits of using axioms as if they were, indeed, freely created from within our own minds. Such habits are now deeply ingrained in the way we teach, as well as in the way we think about and use mathematics. By seldom, if ever, thinking explicitly about where these axioms come from, we can forget the very ground upon which our mathematics rests. Indeed, I recommend re-reading this article multiple times, because, for many, digesting certain key parts of it will be like entering a new town; some things, even important ones, will only be remembered or even noticed on a second or third visit.

---

[26] Specific limit proofs will use this general fact.

[27] Since the real number is defined by the limit, we have to analogically define what we mean by satisfying these axioms using the idea of limit. For instance, each term in the partial sum can be used in the place of the number itself to establish a series for establishing axioms that involve operations such as inversion and multiplication.

[28] Quantity is the first property of physical things; see A. Rizzi, *Physics for Realists: Mechanics* (IAP Press, Baton Rouge, 2008) (see especially chapter I and appendix I).

Yet, even when we forget the foundation, it is no less true that our mathematics stands on the very physical reality from which we pulled it. In fact, our imaginations are constantly at work providing implicit background and foreground pictures that help us think about our mathematics. Though implicit, the principles of reality (with which the axioms of modern mathematics are imbued) are constantly providing the content of what we do, including what consistency any given axiomatic system has.

Modern physics has played a large part in the fecundity of mathematics. It provides a connection to the *full* physical world[29] via analogies to quantity, especially those that robustly link to experiment. However, the success of physics has also been partly due to the axiomatization and symbolization of math and the powerful new way of thinking that comes with it. Indeed, the success of physics could not have happened without it. So, we can say that the success of modern mathematics, including its leaps forward in clarity of proofs (even in those that were already known) and its great discoveries, are in no small part due to the symbolization and the logicization that characterize modern mathematics. Ancient mathematicians did not have this view of mathematics; they kept their methodology (as well as their final aim) tightly tied to the real, which partly explains why mathematics developed so much more slowly in ancient times.

Thus, it would be foolish not to try to get better at teaching students the axioms and the technical skills of modern mathematics. We must, however, not leave students and ourselves without an explicit grounding in what it is all about, or we risk, all of us, missing the whole point of why we do mathematics. Furthermore, we will probably also miss insights into other fields that a closer explicit contact with the source of fecundity of mathematics surely will give. Indeed, we can meet both these goals by more explicitly elucidating the ground from which the axioms spring, because this will make the axiomatic systems and the axioms themselves more accessible to greater numbers and types of people and more understandable to all.

---

[29] Not just the physical world as examined with all properties left aside but extension.

# **Appendix I:** Complex Numbers

In the text, we stopped our process of analogically generalizing the first meaning of number, a whole with parts, after having treated the points on a line as if each was a number. As we have explained, they, of course, are not numbers in this first sense but in an extended sense that is captured in our rule structure visualized and made accessible by our axiom rules. One can go much further with this program of capturing the reality of the first property of physical things into such axiomatic systems. Complex numbers and Quaternions are two obvious examples. We treat here the first, leaving quaternions and other extensions to future publications.

What is left of quantity (extension) to include in our system? The reality of number leaves aside the fact that extended things have parts that share a common boundary. Through the axiomatic formulation discussed in the main text of this article, we imported part of this fuller reality. The number line includes in a certain way, i.e. via a mental construct, the fact that parts lie next to each other in a line. In the world, however, the parts we see are three dimensional. The parts, for example, that make up one's body are three dimensional. Lines have parts in a secondary way. Namely, because they are boundaries of things that have 3-D parts, they, in some way, inherit the divisions of those things.[30] Clearly, to incorporate the fuller reality into our system, we need to somehow incorporate two more dimensions. Complex numbers and quaternions do this.

Complex numbers bring in the second dimension. To do this we have to leave aside the ordering principles introduced in the main body of the article. Our remaining axioms work on a line. Two dimensions can be described by two lines that are perpendicular to one another. To see this, consider moving an object, for example, three units to the right and four units up from a chosen origin. This puts us somewhere in the first quadrant. Physically, this is like considering the fact that any given substance (body) has parts one next to the other and one can count over so many parts and up so many parts to locate a spot in a two dimensional slice of the substance. We can see that finite sized parts will not be enough to locate any spot in the idealized continuum. Said another way, just using integers will not be enough to specify any particular point on the plane. However, it is not hard to see that following the path to irrational numbers to fill in the spots in the line will allow us to locate any point in a plane.

Now, how do we distinguish the horizontal line (*x*-axis) from the vertical? We multiply the vertical number (*y*-axis) by *i* and leave the horizontal number alone. We thus get the standard form of a complex number, $z = a + bi$. This means *a* units along the *x*-axis and *b* units along the *y*-axis.

Next, the question is how to handle the axioms. Start with addition and subtraction, which are straightforward. Start with the location $1+i$ and add $a'+ib'$ to get: $1+a'+i(1+b')$. See how we moved over $a'$ and up $b'$. Subtraction is a similar extension of the ordinary subtraction. The result is the formalism previously developed smoothly takes on this new more general (analogical) meaning for addition and subtraction.

What about multiplication? The following is already clear from the definition of *i* that we have developed so far. Given the point *i* on the *y*-axis, multiplication by "1" gives

---

[30] Bodies have surfaces, i.e. boundaries, that when they are flat and intersect with another surface give a line; consider for example the line that is the place in a rectangular room where two vertical walls meet. In short, a line is the boundary of the boundary of a three dimensional body.

*i.* This means take one *i*, which is one unit along the *y*-axis and no units on the *x*-axis. So we have, after multiplication, a point one unit in the vertical direction along the *y*-axis. So, we have: $1 \times i = i$, which affirms our formal definition of "1" as the identity element. Now, to keep the commutation rule, we want this process to be the same as $i \times 1 = i$, which means: "take 1, *i* times." Given our definition of *i*, this means that "1" on the *x*-axis, gets moved to "1" on the *y*-axis. This can be pictured as follows. A vector from the origin to each point is thus rotated by 90 degrees. That is, the direction from the origin rotates by 90° as shown here for the case of $i \times 1$.

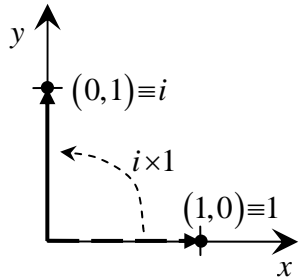

Now, starting with "6", for example, instead of "1" gives $i \times 6$, this again rotates the position of the point by 90°. The commutative rule, $i \times 6 = 6 \times i$, now shows this can also be understood within our formal system as expanding the distance from the origin[31] by a factor of six. Most generally, as we will see in a moment, multiplication consists in a rotation and an expansion (or contraction) of the distance from the origin by a factor. In analogy with the first meaning of multiplication, i.e., a repeated type of addition, we would say $6 \times i$ means "add *i* six times"; the commuted version, $i \times 6$, can mean "add six *i* times" if we take "*i* times" to mean "along the perpendicular axis", so that we have "take six along the *y*-axis." There is a real asymmetry in the way the commuted operations are understood; this is not visible in our formalism, as the multiplication is taken to be commutative.

For now, we have a pressing issue. We know $2 \times 3$ corresponds to no rotation and expanding the length 3 by a factor of 2 to 6. However, what does $i \times i$ mean in our system of axioms? We have seen above that, in the most obvious case $i \times 1$, multiplication by *i* means rotation by 90° but no expansion of the distance. At our level of consideration, there is no fundamental difference between the *y*-axis and the *x*-axis, so hitting *i* with *i* again, should, as before, just rotate without expansion by 90 degrees. This puts us at $-1$. Hence, we take $i \times i = -1$, so that $i \equiv \sqrt{-1}$!

Look what we have done. In order to capture something of the nature of the 2D continuum into our formalism, we have included something we needed to close our formalism in a new way. Up till now, we could define square root for all real numbers, *n*, such that $n \geq 0$; however, there was no number multiplied by itself that could give us negative numbers. The number system was open with respect to the single variable operation of $\sqrt{\phantom{x}}$. This gives us a sign that our formalism is capturing reality in a succinct way, and, even (in an extended sense) a natural way. Similar arguments apply for higher order roots. This leads to the <u>*Fundamental theorem of algebra*</u>: *every degree* $n \geq 1$

---

[31] As illustrated for the previous case in the figure above, this can be represented by a vector when a sense of direction is included.

*polynomial with complex coefficients has exactly n complex roots*.[32] Thus, by incorporating more of extension into our analogical idea of number, we make general solutions to polynomial equations possible.

Again, the inclusion of *i* and the elimination of the order axioms makes our system much more closed while it incorporates new principles of reality from the second dimension. In particular, in our system of axioms the power of a number, i.e.: $x^n$, (which is a special type of multiplication, which, in turn, is a special type of addition) induces us to consider polynomial equations like $x^2 = -1$. This equation is, as we just saw, satisfied by introducing the data from the real world about the second dimension.

Move, now, back to our main path. Look again at the question: what does $i \times i$ mean in terms of the special type of addition analogy? Does it mean, "take *i*, *i* times"? Only weakly. Following our earlier assessment, we take it to mean "take one unit along the axis perpendicular to the perpendicular axis." Again, this is clearly different than simply repeating an addition process. This is the type of analogy that we make to bring further features of 2D into our system, thus making it more general but doing so by analogy. This is why we say it is analogically general.

Next, following our reasoning, multiplying by *i* a third time, so that we have, $i \times i \times i$ should move our point to the negative *y*-axis, to the point: $(0, -1) = 0 - i = -i$. This squares with our formal treatment up to this point. In particular, using $i \times i = -1$ and the associative axiom, $i \times (i \times i) = (i \times i) \times i = -i$. Multiplying a fourth time takes us back to $(1, 0)$.

With this reasoning, we have shown in principle, as promised above, that multiplication corresponds to a rotation followed by an expansion. We can see this explicitly by writing a general complex number in polar form as: $z = r \cos\theta + i\, r \sin\theta$.[33] Clearly, starting with $(1, 0)$, multiplying by *z* takes us an angle $\theta$ in the counterclockwise direction and expands (or contracts) the distance from the origin by *r*.

Division can easily be incorporated as the reverse of multiplication and thus we have all the axioms of a field, so that our new analogical numbers (complex numbers) are fields in which we have introduced the definition of the imaginary number, *i*. We only call *i* "imaginary" because, though it is analogical to the primary meaning of number, it is far removed from it. If properly understood, it is not unreal, just distant from it because of the analogies and axiomatic structure that we create so that it may carry the real aspects of two dimensions that it does within the formalism.

---

[32] This is proved fairly straightforwardly if one starts with the theorem that every order $n \geq 1$ polynomial has at least one prime factor. One first trivially shows first order polynomials have one root and then, using mathematical induction, shows that, given a *k* degree polynomial has *k* roots, then a *k* +1 degree polynomial has *k*+1 roots.

[33] The trigonometric functions involve analogical generalizations as well.

Lastly, using the fact that complex numbers are a field with operations defined above and using the definition of $i$, we recall the Taylor expansion[34] of the exponential to write: $e^{i\theta} = \sum_{n=0}^{\infty} \frac{i^n \theta^n}{n!} = \sum_{n=even}^{\infty} \frac{i^n \theta^n}{n!} + \sum_{n=odd}^{\infty} \frac{i^n \theta^n}{n!} = \sum_{n=0}^{\infty} \frac{i^{2n} \theta^{2n}}{(2n)!} + i\sum_{n=0}^{\infty} \frac{i^{2n} \theta^{2n+1}}{(2n+1)!} =$
$\sum_{n=0}^{\infty} (-1)^n \frac{\theta^{2n}}{(2n)!} + i\sum_{n=0}^{\infty} (-1)^n \frac{\theta^{2n+1}}{(2n+1)!} = \cos\theta + i\sin\theta$. Or, going more directly to the real two dimensional extension, we can add the result term by term using a vector diagram. Take $\theta = 1 \approx \pi/4$, we get: $e^i = 1 + i - \frac{1^2}{2} - \frac{i 1^3}{6} + \frac{1^4}{24} + \frac{i 1^5}{120} - \frac{1^6}{720} - \frac{i 1^7}{5040} + \ldots$. We get:

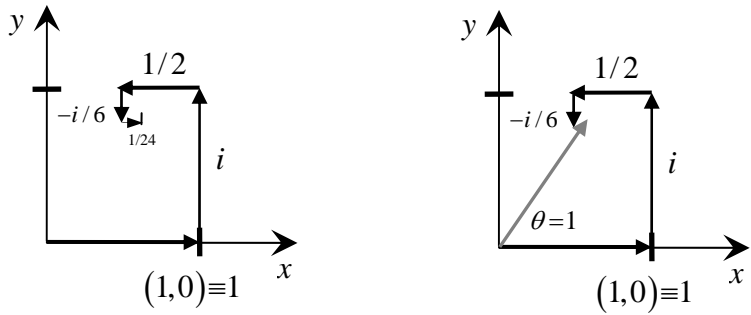

The left figure shows the vectors term by term that we add to get the net result. Note that the vectors wind around in a tighter and tighter circle made from vectors each rotated 90 degrees relative to the one on its tail. The right figure shows the net vector absent some of the vectors shown on the left to avoid the clutter around the arrowhead of the resultant vector, which is shown in gray.

Given this result, we can write any complex in the form $z = re^{i\theta}$. Whether in polar form or "Cartesian" form, the complex number analogically generalizes the concept of number to include important aspects of real extension.

The process of analogical generalization continues in the field of analysis when calculus is applied to real and complex functions.

**References:**
1. Aristotle, edited by J. Barnew, *Complete Works of Aristotle*, Princeton University Press, 1991
2. J. Klein *Greek Mathematical Thought and the Origin of Algebra*, Dover, NY, 1968
3. Dolciani, Beckenbach et al., *Modern Introductory Analysis*, Houghton Mifflin, 1977
4. H.L. Royden, *Real Analysis*, 3rd ed., Stanford University Press, 1968

---

[34] The Taylor expansion and some things that are implicit to it involve analogical generalizations that we have not discussed. We, nonetheless, include it because of the insight in provides to the formalization (encapsulation) of real aspects of extension which are accomplished in modern mathematics.